\documentclass[12pt,a4paper]{article}

\usepackage{graphicx} 
\usepackage[margin=1.2in]{geometry}
\usepackage[font=sf]{caption} 
\usepackage{hyperref}
\usepackage{lipsum}
\usepackage{graphicx}
\usepackage{epstopdf}
\usepackage{algorithmic}
\usepackage{bbm}
\usepackage{amsthm}
\usepackage{amsmath}
\usepackage{amsfonts}
\usepackage{amssymb}
\usepackage{siunitx}
\usepackage{verbatim}
\usepackage{hyperref} 
\usepackage{natbib} 
\usepackage{pifont}
\usepackage{caption}
\usepackage{cleveref}
\newtheorem{theorem}{Theorem}[section]
\newtheorem{corollary}[theorem]{Corollary}
\newtheorem{lemma}[theorem]{Lemma}
\newtheorem{remark}[theorem]{Remark}
\newtheorem{example}[theorem]{Example}
\newtheorem{definition}[theorem]{Definition}
\usepackage{dsfont}
\title{Matrix Concentration: Order versus Anti-order}
\author{ Reihaneh Malekian\thanks{Department of Statistics, Columbia University, New York, NY, USA 
 (\href{mailto:rm3942@columbia.edu}{rm3942@columbia.edu}).}
 \and
 Aaditya Ramdas\thanks{Departments of Statistics and Machine Learning, Carnegie Mellon University Pittsburgh, PA, USA
  (\href{mailto:aramdas@stat.cmu.edu}{aramdas@stat.cmu.edu}).}
}

\begin{document}
\maketitle
\begin{abstract}
   The matrix Markov inequality by \cite{Ahlswede2001} was stated using the Loewner anti-order between positive definite matrices. \cite{Wang2024} use this to derive several other Chebyshev and Chernoff-type inequalities (Hoeffding, Bernstein, empirical Bernstein) in the Loewner anti-order, including self-normalized matrix martingale inequalities. These imply upper tail bounds on the maximum eigenvalue, such as those developed by~\cite{Tropp2012} and~\cite{howard2020time}. The current paper develops analogs of all these inequalities in the Loewner order, rather than anti-order, by deriving a new matrix Markov inequality. These yield upper tail bounds on the \emph{minimum} eigenvalue that are a factor of $d$ tighter than the above bounds on the maximum eigenvalue.
\end{abstract}

\section{Introduction}\label{sec:intro}
  
In this paper, we prove a number of concentration inequalities for random matrices that appear to be new. Inspired by \cite{Ahlswede2001}, we focus on positive semidefinite symmetric random matrices. However, our inequalities are presented using the Loewner order, rather than the anti-order, so our inequalities are neither stronger nor weaker than the original ones. Also, inspired by \cite{Wang2024}, we present randomized versions of some deterministic results which were previously known in the literature, in particular several well known inequalities by \cite{Tropp2012}.

\paragraph{Notation.}
Let $\mathcal{S}_d$ denote the set of all $d \times d$ real-valued symmetric matrices; the set of positive semidefinite and positive definite matrices are denoted by $\mathcal{S}^{+}_d$ and $\mathcal{S}^{++}_d$ respectively. The Loewner partial order is denoted by $ \preceq $, where $A \preceq B$ means $B - A$ is positive semidefinite. We denote the identity matrix by $I$. Consider any function  $f: \mathcal{S}_d \rightarrow \mathcal{S}_d $. 
For a diagonal matrix $\Sigma \in \mathcal{S}_d $, $f (\Sigma)$ is obtained by applying $f$ entrywise; otherwise it is obtained by applying $f$ after orthogonal diagonalization, $f(U\Sigma U) = Uf(\Sigma)U$ where $U$ is orthonormal and $\Sigma$ is diagonal. We shorten \emph{independent and identically distributed} to \emph{iid}. The matrix operator norm of $A$, denoted as $\|A\|$, is its largest singular value.

The central point of the paper will be to compare known inequalities for $\mathbb{P}({X} \npreceq {A})$ to our new ones that are stated for $\mathbb{P}({X} \succeq {A})$. The event $\{ X \npreceq A \}$ means that the largest eigenvalue of $X-A$ is positive, while $\{ X \succeq A \}$ means that the smallest eigenvalue of $X-A$ is nonnegative. Clearly, $\{ X \succeq A \} \subseteq \{ X \npreceq A \}$, and so $\mathbb{P}({X} \succeq {A}) \leq \mathbb{P}({X} \npreceq {A})$. In fact, we will show that the left hand side can be bounded by a factor $d$ times smaller than the corresponding bound for the right hand side.

\section{Main Contributions}

\subsection{Matrix Markov Inequality}

\cite{Ahlswede2001} prove that
$$
\mathbb{P}({X} \npreceq {A}) \leq \operatorname{tr}\left((\mathbb{E} {X}) {A}^{-1}\right),
$$
where ${X}$ is an $\mathcal{S}_d^{+}$-valued random matrix and ${A} \in \mathcal{S}_d^{++}$ is deterministic.

We present the following variant, which replaces the anti-order $\npreceq$ with the order $\succeq$.  Both above and below, it is enough to assume $X \succeq 0$ almost surely.

\begin{theorem}\label{MMI}
Let $X$ be a random matrix with values in $\mathcal{S}_d^{+}$ and let $A \in \mathcal{S}_d^{++}$ be a deterministic matrix. Then,
\begin{equation}\label{eq:matrix-markov2}
\mathbb{P}(X \succeq A) \leq \frac{\operatorname{tr}\left((\mathbb{E} X) A^{-1}\right)}{d} .
\end{equation}
\end{theorem}

The above Theorem~is tight in general, as the following example demonstrates.

\begin{example} Consider a random matrix $X$ and $A \in \mathcal{S}_d^{++}$, such that $X=A$ with probability $p$ and $0$ otherwise. So $\mathbb{E}X= p A$ and $\mathbb{P}(X \succeq A)= \mathbb{P}(X = A)=p$. Also $$\frac{\operatorname{tr}\left((\mathbb{E} X) A^{-1}\right)}{d}= \frac{\operatorname{tr}\left(pA A^{-1}\right)}{d}= \frac{\operatorname{tr}\left(p I\right)}{d}=p.$$
Thus, \eqref{eq:matrix-markov2} holds with equality in this case.
\end{example}

\subsection{Matrix Chebyshev Inequality}
Let $X$ be a random matrix with values in $\mathcal{S}_d$ and $A \in \mathcal{S}_d^{++}$ a deterministic matrix. \cite{Ahlswede2001} showed that
$$\mathbb{P}(|X-\mathbb{E} X| \npreceq A) \leq \operatorname{tr}\left(\left[\mathbb{E}(X-\mathbb{E} X)^2\right] A^{-2}\right).$$ 
This is to be compared with the following theorem.

\begin{theorem}\label{MCI}
Let $X$ be a random matrix with values in $\mathcal{S}_d$ and $A \in \mathcal{S}_d^{++}$ a deterministic matrix. Then,
$$\mathbb{P}(|X-\mathbb{E} X| \succeq A) \leq \frac{1}{d}\operatorname{tr}\left(\mathbb{E}\left[\left(|X-\mathbb{E} X| A^{-1}\right)^2\right]\right)$$
More generally, for $p \geq 1$,
$$\mathbb{P}(|X-\mathbb{E} X| \succeq A) \leq \frac{1}{d}\operatorname{tr}\left(\mathbb{E}\left[\left(|X-\mathbb{E} X| A^{-1}\right)^p\right]\right).$$
Moreover, for $0< q\leq 1$,
$$\mathbb{P}(|X-\mathbb{E} X| \succeq A) \leq \frac{1}{d} \operatorname{tr}\left(\left[\mathbb{E}\left|X-\mathbb{E} X\right|^q\right] A^{-q}\right).$$ 
\end{theorem}

\begin{remark}
$\operatorname{tr}\left(\mathbb{E}\left[\left(|X-\mathbb{E} X| A^{-1}\right)^2\right]\right) = \operatorname{tr}\left(\left[\mathbb{E}(X-\mathbb{E} X)^2\right] A^{-2}\right)$ if $|X-\mathbb{E} X|$ and $A$ commute almost surely.
To see that in general the two quantities are unequal, consider the following example: $X$ is a random matrix that equals $ 
\begin{pmatrix}
-2 & 0 \\
0 & 1
\end{pmatrix}
$ with probability $4/5$ and $ 
\begin{pmatrix}
1 & 0 \\
0 & -10
\end{pmatrix}
$ with probability $1/5$, and $ A = 
\begin{pmatrix}
1 & 0 \\
0 & 1/2
\end{pmatrix}
$.
Then, it is easy to check that $$\operatorname{tr}\left(\mathbb{E}\left[\left(|X-\mathbb{E} X| A^{-1}\right)^2\right]\right) \neq \operatorname{tr}\left(\left[\mathbb{E}(X-\mathbb{E} X)^2\right] A^{-2}\right).$$

\end{remark}

\subsection{Matrix Chernoff Inequalities}
\begin{theorem}\label{MCFI}
Let $X_1\ldots X_n$ be iid matrices with values in $\mathcal{S}_d$, 
and let $A \in \mathcal{S}_d$ and $T, \widetilde{T} \in \mathcal{S}_d^{+}$ be deterministic. Then,
$$\mathbb{P}\left(\sum_{i=1}^n X_i \succeq n A\right) \leq \left\|\mathbb{E}\left(e^{T X \widetilde{T}-T A \widetilde{T}}\right)\right\|^n.$$
\end{theorem}

In contrast, \cite{Ahlswede2001} proved that:
$$
\mathbb{P}\left(\sum_{i=1}^n X_i \npreceq n A\right) \leq d\left\|\mathbb{E}\left(e^{T X \widetilde{T}-T A \widetilde{T}}\right)\right\|^n.
$$


\begin{corollary}\label{corcher}
Let $X_1\ldots X_n$ be iid random matrices, where $0 \preceq X \preceq I$ a.s., and for constants $0 \leq m \leq a \leq1$, we have $\mathbb{E} X \preceq m I$, and $ A \succeq a I$. Then,
$$\mathbb{P}\left(\sum_i^n X_i \succeq n A\right) \leq \exp (-n D(a \| m)),$$ 
where $D(x \| y)=x(\log x-\log y)+(1-x)(\log (1-x)-\log (1-y)).$
\end{corollary}

In contrast, \cite{Ahlswede2001} proved that
$$\mathbb{P}\left(\sum_1^n X_i \npreceq n A\right) \leq d \exp (-n D(a \| m)).$$

\begin{remark}
This Theorem~can be generalized to a similar version of Theorem~5.1 and Corollary 5.2 in \cite{Tropp2012}.
\end{remark}

We end by noting that our results (bounds stated in terms of the order) are dimension-free but anti-order bounds are not.


\subsection{Matrix Laplace Transform Method}
\cite{Tropp2012} proved that
$$
\mathbb{P}\left\{\lambda_{\max }({Y}) \geq t\right\} \leq \inf _{\theta>0}\left\{{e}^{-\theta t} \cdot \mathbb{E} \operatorname{tr} ({e}^{\theta {Y}})\right\},
$$
for all random matrices $Y$ with values in $\mathcal{S}_d$, and all $t \in \mathbb{R}$. In contrast, we show the following two results.

\begin{theorem}\label{laplace}
For any $\mathcal{S}_d$-valued random matrix $Y$ and $t \in \mathbb{R}$, 
 $$
\mathbb{P}\left\{\lambda_{\min }({Y}) \geq t\right\} \leq \frac{1}{d} \inf _{\theta>0}\left\{{e}^{-\theta t} \cdot \mathbb{E} \operatorname{tr} ({e}^{\theta {Y}})\right\}.
$$   
\end{theorem}

\begin{theorem}\label{laplace2}
For any finite sequence $\left\{{X}_k\right\}$ of independent, random matrices with values in $\mathcal{S}_d$ and any $t \in \mathbb{R}$,
$$
\mathbb{P}\left\{\lambda_{\min }\left(\sum_k {X}_k\right) \geq t\right\} \leq \frac{1}{d} \inf _{\theta>0}\left\{{e}^{-\theta t} \cdot \operatorname{tr} \exp \left(\sum_k \log \mathbb{E} {e}^{\theta {X}_k}\right)\right\}.
$$
\end{theorem}

This Theorem~is proved easily using \Cref{laplace} and recalling the following lemma from \cite{Tropp2012}. 

\begin{lemma}
\label{lem:subadditivity-matrix-cgfs}
Consider a finite sequence \(\{X_k\}\) of independent, random, symmetric matrices. Then,
\[
\mathbb{E} \left[ \mathrm{tr} \, \exp \left( \sum_k \theta X_k \right) \right] \leq \mathrm{tr} \, \exp \left( \sum_k \log \mathbb{E} \, e^{\theta X_k} \right) \quad \text{for } \theta \in \mathbb{R}.
\]
\end{lemma}
\begin{remark}
Notice that $\lambda_{\min}(X)= - \lambda_{\max}(-X)$, so the above theorems could be used to derive tail bounds on the spectral norm.
\end{remark}

Also, we can obtain analogous results that correspond to Corollary 3.7, Corollary 3.9, and Theorem~4.1 in \cite{Tropp2012}.

\subsection{Matrix Bernstein Inequalities}

Similar to Theorem~6.1 in \cite{Tropp2012}, we can prove:
\begin{theorem}[Bounded case]\label{Bernstein1}
Consider any finite sequence $\left\{{X}_k\right\}$ of independent, random matrices with values in $\mathcal{S}_d$ such that
$
\mathbb{E} {X}_k={0}$ and $\lambda_{\max }\left({X}_k\right) \leq R$ almost surely.
Then, for all $t \geq 0$,
$$
\begin{aligned}
& \mathbb{P}\left\{\lambda_{\min }\left(\sum_k {X}_k\right) \geq t\right\} \leq  \exp \left(-\frac{\sigma^2}{R^2} \cdot h\left(\frac{R t}{\sigma^2}\right)\right)\leq  \exp \left(\frac{-t^2/2}{\sigma^2+ Rt/3}\right)\\
\end{aligned}
$$
where 
$
\sigma^2:=\left\|\sum_k \mathbb{E}\left({X}_k^2\right)\right\|$ and
$h(x):= (1+x)\log (1+x) - x.$
\end{theorem}
Further, we can show an analogous result for the subexponential case discussed in Theorem~6.2 of \cite{Tropp2012}.
\begin{theorem}[Subexponential case]\label{Bernstein2}
Consider any finite sequence $\left\{{X}_k\right\}$ of independent, random matrices with values in $\mathcal{S}_d$ such that for all $k$ and $p = 2, 3, 4, \ldots$ we have
$$
\mathbb{E} {X}_k={0}, \quad \mathbb{E} (X_k^p) \preceq \frac{p!}{2}\, R^{p-2} A^2_k $$ almost surely.
Then, for all $t \geq 0$,
$$
\begin{aligned}
& \mathbb{P}\left\{\lambda_{\min }\left(\sum_k {X}_k\right) \geq t\right\} \leq  \exp \left(\frac{-t^2/2}{\sigma^2+ Rt}\right),\\
\end{aligned}
$$
where 
$
\sigma^2:=\left\|\sum_k {A}_k^2\right\|$.
\end{theorem}

\subsection{Matrix Azuma Inequality}

\begin{theorem}\label{Azuma}
Consider a finite adapted sequence $\left\{{X}_k\right\}$ of random matrices with values in $\mathcal{S}_d$ and a fixed sequence of $\mathcal{S}_d$-valued matrices $\{A_k \}$, that satisfy
$$
\mathbb{E} ({X}_k \,| \, X_1, \ldots , X_{k-1}) ={0} \quad \text { and } \quad {X^2_k} \preceq A^2_k \quad \text { almost surely. }
$$
Then, denoting $
\sigma^2:=\left\|\sum_k A_k^2 \right\|$, we have that for all $t\geq 0$,
$$
\mathbb{P}\left\{\lambda_{\min} \left(\sum_k X_k \right) \geq t \right\} \leq {e}^{-t^2/{8\sigma^2}},
$$ 
and hence
$$
\mathbb{P}\left\{\lambda_{\max} \left(\sum_k X_k \right) \leq -t \right\} \leq {e}^{-t^2/{8\sigma^2}}.
$$ 
\end{theorem}
In contrast, \cite{Tropp2012} showed that
$$
\mathbb{P}\left\{\lambda_{\max} \left(\sum_k X_k \right) \geq t \right\} \leq d \cdot {e}^{-t^2/{8\sigma^2}},
$$ 
and hence
$$
\mathbb{P}\left\{\lambda_{\min} \left(\sum_k X_k \right) \leq -t \right\} \leq d \cdot {e}^{-t^2/{8\sigma^2}}.
$$ 

\subsection{Matrix McDiarmid Inequality}

\begin{theorem}\label{McDiarmid}
Let $\left\{Z_k: k=1,2, \ldots, n\right\}$ be an independent family of random variables, and let ${H}$ be a function that maps $n$ variables to $\mathcal{S}_d$. Consider a sequence $\left\{{A}_k\right\}$ of fixed matrices in $\mathcal{S}_d$ that satisfy
$$
\left({H}\left(z_1, \ldots, z_k, \ldots, z_n\right)-{H}\left(z_1, \ldots, z_k^{\prime}, \ldots, z_n\right)\right)^2 \preccurlyeq {A}_k^2,
$$
where $z_i$ and $z_i^{\prime}$ range over all possible values of $Z_i$ for each index $i$.
Then, for all $t \geq 0$,
$$
\mathbb{P}\left\{\lambda_{\min }\left({H}(Z)-\mathbb{E} {H}(Z)\right) \geq t\right\} \leq  \mathrm{e}^{-t^2 / 8 \sigma^2}.
$$
where $Z=\left(Z_1, \ldots, Z_n\right)$ and
$\sigma^2:=\left\|\sum_k {A}_k^2\right\|.$
\end{theorem}

In contrast, \cite{Tropp2012} proved that 
$$
\mathbb{P}\left\{\lambda_{\max }\left({H}(Z)-\mathbb{E} {H}(Z)\right) \geq t\right\} \leq d \cdot \mathrm{e}^{-t^2 / 8 \sigma^2}.
$$

\subsection{Matrix Hoeffding's Inequality}

For a random matrix $X$, define the moment-generating function $\psi_X: \mathbb{R}\rightarrow \mathcal{S}_d $ by
$$\psi_X(\lambda) = \mathbb{E}(e^{\lambda X}).$$
From \cite{HDS}, recall that
\begin{definition}[Sub-Gaussian matrices]
\label{def:sub-gaussian-matrices}
A centered symmetric random matrix $X \in \mathcal{S}_{d }$ is sub-Gaussian with parameter $V \in \mathcal{S}^+_{d}$ if, for all $\lambda \in \mathbb{R}$,
$$\psi_X(\lambda) \preceq e^{\frac{\lambda^2 V}{2}}.$$

\end{definition}


\begin{theorem}\label{Hoeffding}
\label{hoeffding}
Let $X_1, \ldots, X_n$ be centered, independent, random matrices in $\mathcal{S}_{d }$ that are sub-Gaussian with parameters $V_1, \ldots, V_n$. Then for all $t \geq 0$, 
we have 
$$
\mathbb{P} \left\{ \lambda_{\min} \left(\frac{1}{n}\sum_{i=1}^n X_i \right) \geq t \right\} \leq \exp \left\{ -\frac{n t^2}{2\sigma^2} \right\},
$$
where $\sigma^2 = \left\| \frac{1}{n}\sum_{i=1}^n V_i \right\|.$
\end{theorem}

In contrast, it is known (restated in \citep{HDS}, Theorem~6.15) that
$$
\mathbb{P} \left\{\left\| \frac{1}{n}\sum_{i=1}^n X_i \right\| \geq t \right\} \leq 2d \exp \left\{ -\frac{n t^2}{2\sigma^2} \right\}.
$$

\subsection{Matrix Maximal Inequalities}

\begin{theorem}[Matrix Doob's Inequality]\label{doob}
Let $\left\{{Y}_n\right\}$ be an $\mathcal{S}_d^{+}$-valued submartingale. For every $N \in \mathbb{N}$ and ${A} \in \mathcal{S}_d^{++}$,
$$
\mathbb{P}\left(\exists n \leq N, {Y}_n \succeq {A}\right) \leq\frac{1}{d} \operatorname{tr}\left(\mathbb{E}\left({Y}_N \mathds{1}_{\left\{\exists n \leq N, {Y}_n \succeq {A}\right\}}\right) {A}^{-1}\right) \leq\frac{1}{d} \operatorname{tr}\left(\left(\mathbb{E} {Y}_N\right) {A}^{-1}\right).
$$
\end{theorem}

In contrast, \cite{Wang2024} proved that,
$$
\mathbb{P}\left(\exists n \leq N, {Y}_n \npreceq {A}\right) \leq \operatorname{tr}\left(\mathbb{E}\left({Y}_N \mathds{1}_{\left\{\exists n \leq N, {Y}_n \npreceq {A}\right\}}\right) {A}^{-1}\right) \leq \operatorname{tr}\left(\left(\mathbb{E} {Y}_N\right) {A}^{-1}\right).
$$

\begin{theorem}[Matrix Ville’s Inequality]\label{Ville}
Let $\{Y_n\}$ be an $\mathcal{S}_d^{+}$-valued supermartingale adapted to some filtration $\{\mathcal{F}_n\}$.  Then, for any $A \in \mathcal{S}_d^{++}$,
$$\mathbb{P}\left(\exists n, Y_{n} \succeq A\right) \leq  \operatorname{tr}((\mathbb{E} Y_0) A^{-1})/d.$$
\end{theorem}
In contrast, \cite{Wang2024} showed that,
$$\mathbb{P}\left(\exists n, Y_{n} \npreceq A\right) \leq  \operatorname{tr}((\mathbb{E} Y_0) A^{-1}).$$

\subsection{Matrix Self-normalized Inequalities}

We can state a lemma similar to the lemma 5.1 in \cite{Wang2024}, as follows:

\begin{lemma}
Let $\left\{{Z}_n\right\}$ be an $\mathcal{S}_d$-valued martingale difference sequence adapted to $\left\{\mathcal{F}_n\right\}$. Let $\left\{{C}_n\right\}$ be an $\mathcal{S}_d$-valued adapted process, $\left\{{C}_n^{\prime}\right\}$ be an $\mathcal{S}_d$-valued predictable process, w.r.t.\ $\left\{\mathcal{F}_n\right\}$. If for all $n$, we have 
$$\mathbb{E}(\exp(Z_n - C_n)\mid \mathcal{F}_{n-1}) \preceq \exp(C^{\prime}_n),$$
then the process
$$
L_n=\operatorname{tr} \exp \left(\sum_{i=1}^n {Z}_i-\sum_{i=1}^n\left(C_i+{C}_i^{\prime}\right)\right)
$$
is a (scalar) supermartingale. Further,
$$
L_n \geq d \cdot \exp \left(\lambda_{\min }\left(\sum_{i=1}^n {Z}_i\right)-\lambda_{\max }\left(\sum_{i=1}^n\left({C}_i+{C}_i^{\prime}\right)\right)\right).
$$
\end{lemma}

\begin{theorem}
[Matrix Empirical Bernstein Inequality]
\label{thm:matrix-eb}
Let $\left\{{X}_n\right\}$ be $\mathcal{S}_d$-valued random matrices adapted to $\left\{\mathcal{F}_n\right\}$ with conditional means $\mathbb{E}\left({X}_n \mid \mathcal{F}_{n-1}\right)={M}_n$. Let $\left\{\widehat{{X}}_n\right\}$ be a sequence of predictable and integrable $\mathcal{S}_d$-valued random matrices such that $\lambda_{\min }\left({X}_n-\widehat{{X}}_n\right) \geq -1$ almost surely. Define $g(x):= -\log{(1-x)}-x$. Then, for any predictable $(0,1)$-valued sequence $\left\{\gamma_n\right\}$,
$$
L_n^{\mathrm{EB}}=\operatorname{tr} \left[\exp \left(\sum_{i=1}^n \gamma_i\left({X}_i-{M}_i\right)-\sum_{i=1}^n g\left(\gamma_i\right)\left({X}_i-\widehat{{X}}_i\right)^2\right)\right]
$$
is a supermartingale. Further, for any stopping time $\tau, \alpha \in(0,1)$, we have
$$
\begin{aligned}\label{eq:matrix-EB}
\mathbb{P}\left(\lambda_{\min }\left(\overline{{X}}_\tau^\gamma-\overline{{M}}_\tau^\gamma\right) \geq \frac{\log ( 1 / \alpha)+\lambda_{\max }\left(\sum_{i=1}^\tau g\left(\gamma_i\right)\left({X}_i-\widehat{{X}}_i\right)^2\right)}{\gamma_1+\cdots+\gamma_\tau}\right) \leq \alpha .
\end{aligned}
$$
\end{theorem}
In contrast, non-randomized version of Theorem~5.6 in \cite{Wang2024} showes us :
$$
\begin{aligned}\label{eq:matrix-EB}
\mathbb{P}\left(\lambda_{\max }\left(\overline{{X}}_\tau^\gamma-\overline{{M}}_\tau^\gamma\right) \geq \frac{\log ( d / \alpha)+\lambda_{\max }\left(\sum_{i=1}^\tau g\left(\gamma_i\right)\left({X}_i-\widehat{{X}}_i\right)^2\right)}{\gamma_1+\cdots+\gamma_\tau}\right) \leq \alpha .
\end{aligned}
$$


\section{Randomized Matrix Concentration Inequalities}

Inspired by \cite{Wang2024}, we now attempt to derive randomized versions of the inequalities in the previous section. In fact, some results are new even without randomization. Recall that a nonnegative random scalar $u$ is  stochastically larger than uniform (\emph{super-uniform}) if $\mathbb{P}(u \leq x) \leq x$ for all $x \geq 0$.


\begin{definition}
We say that a $\mathcal{S}_d^{+}$-valued random matrix ${U}$ is super-uniform if $$\mathbb{P}({U} \preceq Z) \leq \operatorname{tr}(Z)/d$$ for all ${Z} \in \mathcal{S}_d^{+}$.
\end{definition}

If $u$ is a super-uniform scalar, then $u \mathbf{I}$ is a super-uniform matrix, because $\mathbb{P}(u {I} \preceq {Z})=\mathbb{P}\left(u \leq \lambda_{\min }({Z})\right) \leq \lambda_{\min }({Z}) \leq \operatorname{tr}(Z)/d$.
Moreover, we can show for any $Y \in \mathcal{S}_d^{+}$, $U + Y$ is a super-uniform matrix when $U=\sum_i u_i e_i e_i^T$ and $\min_i u_i $ is a super-uniform random scalar.
To see that, consider $Z \in \mathcal{S}_d^{+}$. We can write
$$
\begin{aligned}
& \mathbb{P}(U+Y \preceq Z)\leq \mathbb{P}(U \preceq Z) = \mathbb{P}\left(0\leq \lambda_{\min}(Z-U)\right)\\
& \leq \mathbb{P}\left(\lambda_{\min}(U)\leq \lambda_{\min}(Z)\right) = \mathbb{P}\left(\min_i u_i\leq \lambda_{\min}(Z)\right)\leq \lambda_{\min}(Z)\leq \operatorname{tr}(Z)/d.
\end{aligned}$$

Straightforwardly, the identity matrix ${I}$ is a super-uniform matrix.
In contrast, \cite{Wang2024} call a random matrix $\widetilde{U}$ trace super-uniform if $\mathbb{P}(\widetilde{U} \nsucceq Z) \leq \operatorname{tr}(Z)$ for all ${Z} \in \mathcal{S}_d^{+}$. Notice that this definition is neither stronger nor weaker than our definition of super-uniform matrices and the identity matrix $I$ satistfies both of them.

Now, we state randomized matrix versions of some known theorems. Note that setting $U=I$, as discussed above, recovers the nonrandomized versions.

\begin{theorem}[Randomized Matrix Markov Inequality]\label{RMMI}
Let $X, U$ be two independent random matrices taking values in $\mathcal{S}_d^{+}$, where $U$ is 
super-uniform.  Then, for any $A \in \mathcal{S}_d^{++}$, we have $$\mathbb{P}\left(X \succeq A^{\frac{1}{2}} U A^{\frac{1}{2}}\right) \leq  \operatorname{tr}((\mathbb{E} X) A^{-1})/d.$$
\end{theorem}
In contrast, \cite{Wang2024} showed that for any trace super-uniform  $\widetilde{U}$,
$$\mathbb{P}\left(X \npreceq A^{\frac{1}{2}} \widetilde{U} A^{\frac{1}{2}}\right) \leq  \operatorname{tr}((\mathbb{E} X) A^{-1}).$$
\begin{theorem}[Randomized Matrix Chebyshev Inequality]\label{RMCI}
Let $X, U$ be two independent random matrices taking values in $\mathcal{S}_d$ and $\mathcal{S}_d^{+}$, respectively, where $U$ is 
super-uniform. 
For any $A \in \mathcal{S}_d^{++}$ and $0< q \leq 1$, we have 
$$\mathbb{P}(|X-\mathbb{E} X| \succeq (A^{q/2}UA^{q/2})^{1/q}) \leq \frac{1}{d} \operatorname{tr}\left(\left[\mathbb{E}\left|X-\mathbb{E}X\right|^q\right] A^{-q}\right).$$   
\end{theorem}
A similar result in \cite{Wang2024} shows that for any trace super-uniform matrix $\widetilde{U}$ and $1 \leq p <2$ (note the difference here) we have,

$$\mathbb{P}(|X-\mathbb{E} X| \npreceq(A^{p/2}UA^{p/2})^{1/p}) \leq \operatorname{tr}\left(\left[\mathbb{E}\left|X-\mathbb{E}X\right|^p\right] A^{-p}\right).$$

\begin{theorem}[Randomized Matrix Chernoff Inequality]\label{RMChI}
Let $X$ be a random matrix taking values in $\mathcal{S}_d$ and $U$ a 
super-uniform random matrix taking values in $\mathcal{S}_d^{+}$ independent from $X$. Then, for any $A \in \mathcal{S}_d$ and real number $\gamma$,
$$\mathbb{P}(e^{\gamma X} \succeq e^{ \gamma A/2}Ue^{\gamma A/2 }) \leq \frac{1}{d} \operatorname{tr}\left(\left[\mathbb{E} e^{\gamma X}\right] e^{-\gamma A}\right).$$   
\end{theorem}

\begin{theorem}[Randomized Matrix Chernoff-Hoeffding Inequality]\label{RMCHI}
Let $X_1,\ldots, X_n$ be iid random matrices taking values in $\mathcal{S}_d$ with  mean matrix $\mathbb{E}X_1 = M$, matrix-valued MGF $G(\gamma) = \mathbb{E}e^{\gamma(X_1 - M)}$ and let $U$ be a 
super-uniform random matrix taking values in $\mathcal{S}_d^{+}$ independent from all other matrices. Denoting $\Bar{X}_n := \frac{1}{n}(X_1 +\ldots+X_n)$, we have that for any $\gamma, \beta>0$,
$$\mathbb{P}(e^{\gamma (\Bar{X}_n-M)} \succeq \beta U) \leq \frac{1}{\beta d} \operatorname{tr}\left(e^{n\log{G(\gamma/n)} }\right).$$   
\end{theorem}
Similarly, \cite{Wang2024} showed that for any trace super-uniform $\widetilde{U}$,
$$\mathbb{P}\left(\overline{{X}}_n-{M} \npreceq a{I}+\frac{\log \widetilde{U}}{\gamma}\right) \leq \operatorname{tr}\left({e}^{n \log {G}(\gamma / n)}\right) {e}^{-\gamma a}.$$
 However, the two results are not directly comparable because the matrix logarithm is operator monotone, meaning that $\mathbb{P}(e^{\gamma (\Bar{X}_n-M)} \succeq \beta U) \leq \mathbb{P}(\gamma (\Bar{X}_n-M) \succeq \log (\beta U) )$, which may then not be upper bounded by the same right hand side. Nevertheless, the comparison to their result becomes easier to intuit by setting $a = \log \beta / \gamma$.

\begin{theorem}[Randomized Matrix Ville’s Inequality]\label{RVille}
Let $\{Y_n\}$ be an $\mathcal{S}_d^{+}$-valued supermartingale adapted to some filtration $\{\mathcal{F}_n\}$, and let $\tau$ be a stopping time on $\{\mathcal{F}_n\}$. Let $U$ be a 
super-uniform random matrix taking values in $\mathcal{S}_d^{+}$ independent from $\mathcal{F}_{\infty}$. Then, for any $A \in \mathcal{S}_d^{++}$,
$$\mathbb{P}\left(Y_{\tau} \succeq A^{\frac{1}{2}} U A^{\frac{1}{2}}\right) \leq  \operatorname{tr}((\mathbb{E} Y_0) A^{-1})/d.$$
\end{theorem}
Theorem~\ref{Ville} follows as a consequence of the above. In contrast, \cite{Wang2024} showed that for any trace super-uniform matrix $\widetilde{U}$,
$$\mathbb{P}\left(Y_{\tau} \npreceq A^{\frac{1}{2}} \widetilde{U} A^{\frac{1}{2}}\right) \leq  \operatorname{tr}((\mathbb{E} Y_0) A^{-1}).$$

\begin{remark}
We can prove several more results analogous to what has been shown in all aforementioned papers and also the ones in \cite{Aoun2019}. 
\end{remark}

\begin{theorem}
[Uniformly Randomized Matrix Empirical Bernstein Inequality]
Under the same setting as Theorem~\ref{thm:matrix-eb}, 
for any super-uniform scalar random variable $U$ independent of $\mathcal{F}_{\infty}$,~\eqref{eq:matrix-EB} holds with $\log(U/\alpha)$ in place of $\log(1/\alpha)$.
\end{theorem}
\begin{remark}
We can also prove a result similar to Theorem~5.5 in \cite{Wang2024}.
\end{remark}

\section{Proofs}
\subsection{Proof of \Cref{MMI}}
We provide the proof here for the case that $X$ has a discrete distribution for transparency;  the same argument could be rewritten for more general cases.

\begin{proof}
Define $Y=A^{-\frac{1}{2}} X A^{-\frac{1}{2}}$. Clearly $\{X \succeq A\}=\{Y  \succeq I\}$.
Now,
$$
\begin{aligned}
& A^{-\frac{1}{2}}(\mathbb{E} X) A^{-\frac{1}{2}}=\mathbb{E} Y= \sum_y \mathbb{P}(Y=y) y \succeq \sum_{y \succeq I} \mathbb{P}(Y=y) y \succeq \sum_{y \succeq I} \mathbb{P}(Y=y) I =  \mathbb{P}(Y \succeq I) I. 
\end{aligned}
$$
By trace rotation and monotonicity, we obtain 
\[ \operatorname{tr}\left((\mathbb{E} X) A^{-1}\right)=\operatorname{tr}\left(A^{-\frac{1}{2}}(\mathbb{E} X) A^{-\frac{1}{2}}\right) = \operatorname{tr}(\mathbb{E} Y) \geq d \cdot \mathbb{P}(Y \succeq I),
\]
proving the claim.
\end{proof}

\subsection{Proof of \Cref{MCI}}
\begin{proof}
Write
$$\mathbb{P}(|X-\mathbb{E} X| \succeq A)= \mathbb{P}(|X-\mathbb{E} X|A^{-1} \succeq I)\leq  \mathbb{P}(\left(|X-\mathbb{E} X|A^{-1}\right)^{p} \succeq I),$$
and then apply the Matrix Markov Inequality.
For $0< q\leq 1$, notice that $X \mapsto X^q$ is a monotone function with respect to $\succeq$ and then apply the Matrix Markov Inequality.
\end{proof}

\subsection{Proof of \Cref{MCFI}}
\begin{proof}

First notice that for any random matrix $Y \in \mathcal{S}_d$ and any deterministic matrix $B \in \mathcal{S}_d$, we can write
$$
\begin{aligned}
& \mathbb{P}(Y \succeq B) = \mathbb{P}(Y-B \succeq 0) = \mathbb{P}\left(T(Y-B)\widetilde{T} \succeq 0\right) \leq \mathbb{P}\left( e^{T\left(Y-B\right)\widetilde{T}} \succeq I \right).
\end{aligned}$$
Applying \Cref{MMI} gives us a further upper bound of
$$\mathbb{P}(Y \succeq B) \leq \frac{1}{d}\mathrm{tr} \left( \mathbb{E} \left( e^{T Y \widetilde{T}- T B \widetilde{T}}\right) \right).$$
So,
\[
\mathbb{P}\left( \sum_{i=1}^n X_i \succeq nA \right) \leq \frac{1}{d}\mathrm{tr} \left( \mathbb{E} \left( e^{T \sum_{i=1}^n X_i \widetilde{T}- T n A \widetilde{T}}\right) \right).
\]
Now, recall the Golden-Thompson inequality \citep{Petz1994}: For \(A, B \in \mathcal{S}_d\),
\(
\mathrm{tr} \left( e^{A + B} \right) \leq \mathrm{tr} \left( e^A e^B \right).
\)
Thus, we can iteratively simplify the above upper bound as
\[
\begin{aligned}
& \frac{1}{d}\mathrm{tr} \left( \mathbb{E} \left( e^{T \sum_{i=1}^n X_i \widetilde{T}- T n A \widetilde{T}}\right) \right) = \frac{1}{d} \mathbb{E} \left( \mathrm{tr} \left( e^{T \sum_{i=1}^n X_i \widetilde{T}- T n A \widetilde{T}}\right)  \right)\\
& \leq \frac{1}{d} \mathbb{E} \left( \mathrm{tr} \left( e^{T \sum_{i=1}^{n-1} X_i \widetilde{T}- T (n-1) A \widetilde{T}}  e^{T(X_{n} - A)\widetilde{T}} \right)\right) \\
& = \frac{1}{d}\mathbb{E}\left( \mathbb{E} \left[ \mathrm{tr} \left( e^{T \sum_{i=1}^{n-1} X_i \widetilde{T}- T (n-1) A \widetilde{T}}  e^{T(X_{n} - A)\widetilde{T}} \right)\mid X_1, \ldots, X_{n-1}\right]\right)  \\
& = \frac{1}{d}\mathbb{E}\left(  \mathrm{tr}\left[ \mathbb{E} \left( e^{T \sum_{i=1}^{n-1} X_i \widetilde{T}- T (n-1) A \widetilde{T}}  e^{T(X_{n} - A)\widetilde{T}} \right)\mid X_1, \ldots, X_{n-1}\right]\right)  \\
& = \frac{1}{d}\mathbb{E}\left(  \mathrm{tr}\left[ e^{T \sum_{i=1}^{n-1} X_i \widetilde{T}- T (n-1) A \widetilde{T}} \mathbb{E} \left(   e^{T(X_{n} - A)\widetilde{T}} \right)\right]\right)  \\
& \leq \frac{1}{d} \left\| \mathbb{E} \left(e^{T(X_{n} - A)\widetilde{T}}\right) \right\| \mathbb{E} \left( \mathrm{tr} \left( e^{T \sum_{i=1}^{n-1} X_i \widetilde{T}- T (n-1) A \widetilde{T}}\right) \right) \\
& \leq \ldots \leq \frac{1}{d}\left\| \mathbb{E} \left( e^{T(X- A)\widetilde{T}} \right) \right\|^n \mathbb{E}\left(\mathrm{tr}(I)\right)= \left\| \mathbb{E} \left( e^{T(X- A)\widetilde{T}} \right) \right\|^n .
\end{aligned}
\]

\end{proof}

\subsection{Proof of \Cref{corcher}}
\begin{proof}
Apply \Cref{MCFI} with $T = \widetilde{T}= \sqrt{t} I$, the rest is similar to the proof of Theorem 19 in \cite{Ahlswede2001}.
\end{proof}

\subsection{Proof of \Cref{laplace}}
\begin{proof}
For any $\theta > 0$,
$$\mathbb{P}\left\{\lambda_{\min }({Y}) \geq t\right\} = \mathbb{P}\left\{e^{\theta \lambda_{\min }({Y})} \geq e^{\theta t}\right\} \leq e^{-\theta t} \mathbb{E}(e^{\theta \lambda_{\min }({Y})}),$$
where we used the usual Markov inequality.
Now notice that $$e^{\theta \lambda_{\min }({Y})}= \lambda_{\min }(e^{\theta Y}) \leq \frac{1}{d} \mathrm{tr}(e^{\theta Y}),$$ so
$$ e^{-\theta t} \mathbb{E}(e^{\theta \lambda_{\min }({Y})}) \leq \frac{e^{-\theta t}}{d} \mathbb{E}\left(\mathrm{tr}(e^{\theta Y})\right). $$
\end{proof}

\subsection{Proof of \Cref{Bernstein1}}
\begin{proof}
We skip this proof, since it is very similar to the proof of Theorem~6.1 in \cite{Tropp2012}.
\end{proof}

\subsection{Proof of \Cref{Bernstein2}}
\begin{proof}
We skip this proof, since it is very similar to the proof of Theorem~6.2 in \cite{Tropp2012}.
\end{proof}

\subsection{Proof of \Cref{Azuma}}
\begin{proof}
The proof of the first inequality is similar to the proof of Theorem~7.1 in \cite{Tropp2012}, applying \cref{laplace2}.
To get the second inequality, we apply the Theorem~on $\{- X_k \}$.
\end{proof}
\subsection{Proof of \Cref{McDiarmid}}
This proof is similar to the proof of the corollary 7.5 in \cite{Tropp2012}, we need to use \cref{Azuma} to complete it.

\subsection{Proof of \Cref{hoeffding}}
\begin{proof}
Since $X \mapsto \log(X)$ is monotone, for all $i$ and constants $a_i$ we have
$$
\begin{aligned}
& \log \psi_{X_i}\left(a_i \theta\right)\preceq  \frac{a_i^2\theta^2 }{2}V_i \\
& \Rightarrow  \sum_i \log \psi_{X_i}\left(a_i \theta\right)\preceq \frac{\theta^2 }{2}\sum_i a_i^2V_i \\
& \Rightarrow \operatorname{tr}\left(\exp\left\{\sum_i \log \psi_{X_i}\left(a_i \theta\right)\right\}\right) \leq \operatorname{tr}\left(\exp\left\{\frac{\theta^2 }{2}\sum_i a_i^2V_i\right\}\right). \\
\end{aligned}
$$
The last inequity is justified by trace monotoniciy, which is as follows:

Let $f: \mathbb{R} \rightarrow \mathbb{R}$ be an increasing function. Then so is $\operatorname{tr} \circ f: \mathcal{S}_d \rightarrow \mathbb{R}$, i.e.,
$$
{A} \preceq {B} \Longrightarrow \operatorname{tr} f({A}) \leq \operatorname{tr} f({B}).
$$

Now, using \Cref{laplace2},
$$
\begin{aligned}
& \mathbb{P}\left(\lambda_{\min }\left(\sum a_i X_i\right) \geq t\right) \leq \frac{1}{d} \inf _{\theta>0}\left\{e^{-\theta t}  \operatorname{tr}\left(\exp\left\{\sum_i \log \psi_{X_i}\left(a_i \theta\right)\right\}\right)\right\} \\
& \leq \frac{1}{d} \inf _{\theta>0}\left\{e^{-\theta t} \operatorname{tr}\left(\exp\left\{\frac{\theta^2 }{2}\sum_i a_i^2V_i\right\}\right)\right\} \\
& \leq \frac{1}{d} \inf _{\theta>0}\left\{e^{-\theta t} d\, \lambda_{\max }\left(e^{\frac{\theta^2}{2} \sum a_i^2 V_i}\right)\right\} \\
& =\inf_{\theta>0}\left\{e^{-\theta t} e^{\frac{\theta^2}{2}\left\| \sum_{i=1}^n a_i V_i \right\| }\right\}.
\end{aligned}
$$
If $\theta={t}/{\left\| \sum_{i=1}^n a_i V_i \right\|}$, and for all $i$, $a_i = 1/n$, we get the desired inequality.
\end{proof}
\subsection{Proof of \Cref{doob}}
\begin{proof}
It is similar to the proof of Theorem~4.14, in \cite{Wang2024}, so we skip.

\end{proof}
\subsection{Proof of \Cref{RMMI}}
\begin{proof}
By direct calculation, we have
    \[
\mathbb{P} \left( X \succeq A^{1/2} U A^{1/2} \right) 
= \mathbb{P} \left( A^{-1/2} X A^{-1/2} \succeq U \right)
= \mathbb{E} \left( \mathbb{P} \left( A^{-1/2} X A^{-1/2} \succeq U \mid X \right) \right)
\]
\[
\leq \frac{1}{d} \mathbb{E} \left( \mathrm{tr} \left( A^{-1/2} X A^{-1/2} \right) \right) = \frac{1}{d}\mathrm{tr} \left( (\mathbb{E} X) A^{-1} \right),
\]

\end{proof}
\subsection{Proof of \Cref{RMCI}}
Apply \Cref{RMMI} and notice that $X \mapsto X^q$ is monotone for $0<q \leq 1$.

\subsection{Proof of \Cref{RMCHI}}
By \Cref{RMChI}, we have 
$$\mathbb{P}(e^{\gamma (\Bar{X}_n-M)} \succeq e^{ \gamma A/2}Ue^{\gamma A/2 }) \leq \frac{1}{d} \operatorname{tr}\left(\left[\mathbb{E} e^{\gamma (\Bar{X}_n-M)}\right] e^{-\gamma A}\right). $$

If $A := \log(\beta I)/{\gamma},$

$$
\begin{aligned}
& \mathbb{P}(e^{\gamma (\Bar{X}_n-M)} \succeq \beta U) \leq \frac{1}{\beta d} \operatorname{tr}\left(\left[\mathbb{E} e^{\gamma (\Bar{X}_n-M)}\right]\right)=\frac{1}{\beta d} \mathbb{E}\left(\left[ \operatorname{tr}(e^{\gamma (\Bar{X}_n-M)})\right]\right)\\
& =\frac{1}{\beta d} \mathbb{E}\left(\left[ \operatorname{tr}(e^{\frac{\gamma}{n} (\sum_{i=1} ^n {X_i}-nM)})\right]\right) \leq \frac{1}{\beta d} \mathbb{E}\left(\left[ \operatorname{tr}(\prod_{i=1}^n e^{\frac{\gamma}{n} ({X_i}-M)})\right]\right)=\frac{1}{\beta d} \operatorname{tr}\left[\mathbb{E} \prod_{i=1}^n e^{\frac{\gamma}{n} ({X_i}-M)}\right]\\
& =\frac{1}{\beta d} \operatorname{tr}\left[ \prod_{i=1}^n \mathbb{E} e^{\frac{\gamma}{n} ({X_i}-M)}\right] =\frac{1}{\beta d} \operatorname{tr}\left(e^{n \log{G(\gamma/{n}})}\right).
\end{aligned}
$$

\subsection{Proof of \Cref{RVille}}
It is similar to the proof of Theorem~4.4, in \cite{Wang2024}, so we skip.

\subsection*{Acknowledgments} 
We thank Joel Tropp and Hongjian Wang for helpful discussions.

\bibliographystyle{apalike}
\bibliography{example}

\end{document}